\documentclass[final,reqno]{article}
\usepackage{amsmath, amssymb}
\usepackage{amsthm}
\usepackage[table]{xcolor}
\usepackage{subfig}
\usepackage{graphicx}
\usepackage{epstopdf}
\usepackage{amsfonts}
\usepackage{tikz}
\usepackage[table]{xcolor}
\usepackage{booktabs}
\usepackage{epsfig,epsf,psfrag}
\usepackage[left=1in, right=1in, top=1in, bottom=1in]{geometry}
\newtheorem{remark}{Remark}

\newtheorem{cor}{Corollary}
\newtheorem{lemma}{Lemma}
\newtheorem{theorem}{Theorem}
\DeclareMathOperator{\sgn}{sgn}
\DeclareMathOperator{\erfc}{erfc}

\newcommand{\myover}[2]{\genfrac{}{}{0pt}{}{#1}{#2}}
\newcommand{\be} [1] {\begin{equation} \label{#1}}
\newcommand{\ee}{\end{equation}}
\def\x{{\bf x}}
\def\y{{\bf y}}

\def\dt{{\Delta t}}

\def\n{{\bf n}}

\def\qed{{\hfill \Box}}
\def\beq{\begin{equation}}
\def\eeq{\end{equation}}

\def\bR{{\bf R}}

\def\bN{{\bf N}}

\def\rprime{\hbox{\hskip.25em\raise.5ex\hbox{$'$}\hskip.15em}}

\def\ddn1{{\frac{\partial}{\partial \nu_{\yb}}}}

\title{Hybrid asymptotic/numerical methods for the evaluation of 
layer heat potentials in two dimensions}
\author{Jun Wang%
\thanks{Courant Institute of Mathematical Sciences,
New York University, New York, New York 10012.
Present address: Flatiron Institute,
Simons Foundation, New York, New York 10010.
Email: jwang@flatironinstitute.org.}
\and
Leslie Greengard%
\thanks{Courant Institute of Mathematical Sciences,
New York University, New York, New York 10012 and
Flatiron Institute, Simons Foundation, New York, New York 10010.
This work was supported in part by the Applied Mathematical
Sciences Program of the U.S. Department of Energy under Contract
DEFGO288ER25053 and by the RiskEcon Lab for Decision Metrics
at the Courant Institute.
Email: greengard@cims.nyu.edu, Tel: (212) 998-3306,
Fax: (212) 995-4121.}
}

\begin{document}

\maketitle
\pagestyle{myheadings}
\markboth{\sc J. Wang and L. Greengard}
{\sc Asymptotic/numerical methods for layer heat potentials}

\begin{abstract}
We present a hybrid asymptotic/numerical method for the accurate 
computation of single and double layer heat
potentials in two dimensions. It has been shown in previous work that
simple quadrature schemes 
suffer from a phenomenon called 
``geometrically-induced stiffness,"  meaning that formally high-order
accurate methods require excessively small time steps before the 
rapid convergence rate is observed.
This can be overcome by analytic integration in time,
requiring the evaluation of a collection of spatial boundary 
integral operators with non-physical, weakly singular kernels. 
In our hybrid scheme, we combine a local asymptotic approximation with
the evaluation of a few boundary integral operators involving only
Gaussian kernels, which 
are easily accelerated by a new version of the fast Gauss transform.
This new scheme is robust, avoids
geometrically-induced stiffness, and is 
easy to use in the presence of moving geometries. Its 
extension to three dimensions is natural and straightforward, and
should permit layer heat potentials to become flexible and powerful 
tools for modeling diffusion processes. 
\end{abstract}


\vspace{.1in}



\section{Introduction}

A variety of problems in computational physics involve the 
solution of the heat equation in moving geometries - often
as part of a more complex modeling task. Examples include 
heat transfer, fluid dynamics, solid-liquid phase transitions, and the
diffusion of biochemical species.
Integral equation methods are particularly powerful for solving such problems;
they are stable, insensitive to the complexity of the geometry,
and naturally applicable to non-stationary domains.
For unbounded domains, they have the additional advantage that they do not
require artificial absorbing conditions to be imposed on a finite
computational domain, as do methods based on directly discretizing the 
governing partial differential equation.
Methods of this type have not become widespread, in part because of 
quadrature difficulties when evaluating the necessary space-time
integrals. 
Here, we present a robust, new method that overcomes these difficulties, 
combining a local asymptotic calculation with an exponential 
variable transformation that permits high order accuracy to be achieved 
using a small number of quadrature nodes in time, 
{\em both on and off surface}.

Before turning to the more technical aspects of quadrature, let us
introduce the issues with a simple model problem.
We seek the solution of the heat equation 
in either the interior or exterior of 
a non-stationary domain given at time $t$ by $\Omega(t)$.
The boundary of $\Omega(t)$ will be denoted by 
$\Gamma(t)$. More precisely, we seek to compute the solution to
the system:

\begin{alignat}{4}
\frac{\partial U}{\partial t}(\x, t)-\nabla^2 U(\x,t)& = 0, 
\quad & (\x,t) & \in \Omega(t)
\label{heatfree} \\
U(\x,0) & = f(\x), \quad & \x & \in \Omega(0),
\label{eq:init} \\
 \left[ U(\x,t) \right] & = \mu(\x,t), 
\quad & (\x,t) & \in \Gamma(t)
\label{eq:bc1} \\
 \left[ \frac{\partial U}{\partial n}(\x,t) \right] & = \sigma(\x,t), 
\quad & (\x,t) & \in \Gamma(t)
\label{eq:bc2}
\end{alignat}
over the time interval $0 \leq t \leq T$.
Here, the notation $\left[ f(\x,t) \right]$ denotes the jump in
$f(\x,t)$ across $\Gamma(t)$ in the outward normal direction.
We let $N$ denote the total number of time steps 
with $\dt = T/N$.
Without loss of generality, let us simply consider the first time step.
Given the solution to \eqref{heatfree} at time $t_0=0$,
the solution at time $t=t_0+\dt = \dt$ can be expressed 
explicitly as
\begin{equation}
u(\x,\dt) = J[U(\x,0),\dt](\x) + D[\mu](\x,\dt) + S[\sigma](\x,\dt), 
\label{fullsol}
\end{equation}
with 
\begin{eqnarray}
J[f,\dt](\x) &=& \int_{\Omega(0)} G(\x-\y,\dt) f(\x) \, d\y \, ,
\label{potdef} \\
S[\sigma](\x,\dt) &=& 
\int_{0}^{\dt}\int_{\Gamma(\tau)} G(\x-\y,\dt-\tau) 
\sigma(\y,\tau) ds_{\y}d\tau \, ,
\label{slpdef} \\
D[\mu](\x,\dt) &=& \int_{0}^{\dt}\int_{\Gamma(\tau)} 
\frac{\partial}{\partial \nu_{\y}} G(\x-\y,\dt-\tau)
\mu(\y,\tau) ds_{\y}d\tau \, .
\label{dlpdef} \\
\end{eqnarray}
Here,
\begin{equation}
G(\x,t) = \frac{e^{-\|\x\|^2/4t}}{(4 \pi t)^{d/2}} 
\label{heatker}
\end{equation}
is the fundamental solution of the heat equation in $d$ dimensions.
The function $J[f,\dt](\x)$ is referred to as 
an {\em initial} (heat) potential,
the function $S[\sigma](\x,\dt)$ is referred to as a {\em single layer}
potential, and
the function $D[\mu](\x,\dt)$ is referred to as a {\em double layer}
potential.
In the remainder of this paper, we will assume $d=2$.

Both $S[\sigma](\x,\dt)$ and
$D[\mu](\x,\dt)$ satisfy the homogeneous heat equation.
The fact that \eqref{fullsol} is the exact solution to the model
problem is a consequence of the following theorem concerning the 
behavior of layer potentials 
\cite{guentherlee,pogorzelski}.

\begin{theorem}
Let $\Omega(t)$ be a bounded domain with smooth 
boundary $\Gamma(t)$, and let $\x'$ denote a point on $\Gamma(t)$.
Then 
$D[\mu](\x,0) = 0$ for $\x \notin \Gamma(0)$, and 
\begin{eqnarray}
\lim_{\myover{\x \rightarrow \x'}{\x \in \Omega(\dt)}}
D[\mu](\x,\dt) &=& 
-\frac{1}{2}\mu(\x',\dt) + D^{*}[\mu](\x',\dt) \, , \label{dlpin} \\
\lim_{\myover{\x \rightarrow \x'}{\x \in c\Omega(\dt)}}
D[\mu](\x,\dt) &=& 
\frac{1}{2}\mu(\x',\dt) + D^{*}[\mu](\x',\dt) \, ,  \label{dlpout}
\end{eqnarray}
where 
\be{eq:dl}
D^{*}[\mu](\x',\dt) := 
\int \limits_{0}^\dt \int \limits_{\Gamma(\tau)}\frac{\partial}{\partial n_{\y}}
G(\x'-\y,\dt-\tau) \; \mu(\y,\tau) \; ds_{\y} \;d\tau, \quad \x' \in \Gamma(\dt)
\ee
is a weakly singular operator acting on the boundary, and
$c\Omega(\dt)$ denotes the complement of $\Omega(\dt)$.

We also have
$S[\sigma](\x,0) = 0$ for $\x \notin \Gamma(0)$, and 
\begin{eqnarray}
\lim_{\myover{\x \rightarrow \x'}{\x \in \Omega(\dt)}}
\frac{\partial}{\partial n_{\x'}} S[\sigma](\x,\dt) &=& 
\frac{1}{2}\sigma(\x',\dt) + K^{*}[\sigma](\x',\dt) \, , \label{slpin} \\
\lim_{\myover{\x \rightarrow \x'}{\x \in c\Omega(\dt)}}
\frac{\partial}{\partial n_{\x'}} S[\sigma](\x,\dt) &=& 
-\frac{1}{2}\sigma(\x',\dt) + K^{*}[\sigma](\x',\dt) \, ,  \label{slpout}
\end{eqnarray}
where 
\be{eq:slnormal}
K^{*}(\sigma)(\x',\dt) := 
\int \limits_{0}^\dt \int \limits_{\Gamma(\tau)}
\frac{\partial}{\partial n_{\x'}}
G(\x'-\y,\dt-\tau) \; \sigma(\y,\tau) \; ds_{\y} \;d\tau, \quad \x' \in \Gamma(\dt),
\ee
is a weakly singular operator acting on the boundary.

Finally, the kernel of the single layer
potential is also weakly singular, so that
\be{eq:sl}
S[\sigma](\x',\dt) := 
\int \limits_{0}^\dt \int \limits_{\Gamma(\tau)} G(\x'-\y,\dt-\tau) \; 
\sigma(\y,\tau) \; ds_{\y} \;d\tau
\ee
is well-defined for $\x' \in \Gamma(\dt)$, with 
$\lim_{\x \rightarrow \x'} S[\sigma](\x,\dt) = 
S[\sigma](\x',\dt)$.
\end{theorem}

\begin{remark}
The evaluation of $J[f,\dt]$ can be carried out using a ``volume
integral" version of the fast Gauss transform \cite{fgtvol}, so we
concentrate here on the efficient and accurate evaluation of the
boundary integral components
$D[\mu](\x,\dt)$ and $S[\sigma](\x,\dt)$.
\end{remark}

\begin{remark}
Boundary integral methods that
use \eqref{fullsol} over  the entire space-time history of the problem
without recasting the solution 
as new initial value problem at each time step are also very powerful,
but require fast algorithms that
are outside the scope of the present paper 
\cite{greengard_lin,greengard_strain2,tauschheat1}. Even in those schemes, 
the accurate computation of $S[\sigma](\x,\dt)$ and 
$D[\mu](\x,\dt)$ is essential.
\end{remark}

\begin{remark}
For Dirichlet, Neumann, or Robin boundary conditions on $\Gamma(t)$,
classical approaches based on the single and/or double layer potential
yield Volterra integral equations of the second kind 
\cite{guentherlee,pogorzelski}. We omit a discussion of these boundary
value problems here, since we are focused on the problem of quadrature,
and the relevant considerations are more easily illustrated by considering
the analytic solution
\eqref{fullsol} with $\sigma$ and $\mu$ viewed as known functions.
\end{remark}

\section{Asymptotics of local heat potentials in two dimensions} 
\label{sec_asymp}

We turn first to the construction of
asymptotic expansions for $S[\sigma]$ and $D[\mu]$ that are 
valid both on and off
surface.  Earlier results
in \cite{greengard_strain2} can be obtained as the limiting case for $S[\sigma]$ when  
the target approaches the boundary.
Similar analysis (for the single layer case) was carried
out by Strain \cite{strain_adapheat}, using a different set of asymptotic
parameters. 

Let us assume that an arbitrary target point $\x$ is expressed
in the form $\x=\x_0+r\cdot\n$, 
where $\x_0\in\Gamma(\dt)$ is the closest point on the boundary at time 
$\dt$, $\n$ is the unit inward normal vector
at $\x_0$, and $r$ is the signed distance from $\x$ to the boundary.
That is, when $\x$ is in the interior of the domain, $r$ is taken to 
be positive, while if $\x$ is an exterior point, then $r$ is negative.

Observing that the single and double layer potentials are invariant
under Euclidean motion, we shift and rotate 
the coordinate system so that $\x_0$ lies at the
origin $(0,0)$ and $\n$ points in direction of the positive $y$-axis. 
We further assume that the boundary $\Gamma(\dt)$ extends 
away from the origin in such a way that
it is locally the graph of a function, parametrized as
\begin{equation}
x=s,\;\;\;y=y(s,\tau)
\label{cofv}
\end{equation}
where $y(0,\dt)=y_s(0,\dt)=0$ with curvature $\kappa=y_{ss}(0,\dt)$ 
and normal velocity $v=y_{\tau}(0,\dt)$.

\begin{remark}
With a slight abuse of notation, we will assume that the boundary 
extends to infinity while remaining the graph of a function.
This involves an exponentially small error, assuming that
$\Gamma(\dt)$ is the graph of a function out to a distance 
$L$, where $L > \sqrt{\dt} \ln(1/\epsilon)$. From the decay
of the heat kernel, it is straightfoward to verify that this incurs
an error of the order $O(\epsilon)$, so long as the boundary then extends 
farther away from the target point $\x_0$, now located at the origin.
\end{remark}

For the sake of clarity, let us now write $\sigma(s,\tau)$ in place of 
$\sigma(\y,\tau)$ and switch to the coordinate system indicated by
\eqref{cofv}.
The single layer potential then takes the form:
\begin{equation}
S[\sigma](\x,\dt)\approx \int_{0}^\dt \int_{-\infty}^{\infty}
\frac{e^{-s^2/4(\dt-\tau)}e^{-(r-y)^2/4(\dt-\tau)}}{4\pi(\dt-\tau)}\sigma(s,\tau)
\sqrt{1+y_s^2(s,\tau)} \, ds \, d\tau,
\end{equation}
where $y = y(s,\tau)$.
If we apply the change of variables
\begin{equation}
z=\sqrt{4(\dt-\tau)},\;\;\;u=\frac{s}{\sqrt{4(\dt-\tau)}} \, ,
\end{equation}
we obtain
\begin{equation}\label{eq:simplifiedsloc}
S[\sigma](\x,\dt)=\frac{1}{2\pi}\int_0^{2\sqrt{\dt}}\int_{-\infty}^{\infty}e^{-u^2}e^{-(r-y)^2/z^2}\sigma(uz,\dt-z^2/4)\sqrt{1+y_s^2} \, du \, dz .
\end{equation}
From this formula, 
we expand $y(s,\tau)$ and $\sigma(uz,\dt-z^2/4)$ as Taylor series in 
$u$ and $z$.
This yields asymptotic expansions for the local parts of layer potentials. 

\begin{lemma}\label{lemma:slocasym}
Let $\Gamma(\tau)$ and $\sigma(\x,\tau)$ be four times differentiable. Then
\begin{equation}\label{eq:slocasym}
S[\sigma](\x,\dt)=\frac{1}{2}
\sqrt{\frac{\dt}{\pi}}E_{3/2}\left(\frac{c^2}{4}\right)
\left(1+\frac{\kappa-v}{2}\cdot c\sqrt{\dt}\right)
\cdot\sigma(\x_0,\dt)+O(\dt^{3/2})
\end{equation}
where $\x=\x_0+c\sqrt{\dt}\cdot\n$, $\x_0\in\Gamma(\dt)$ 
is a point on the boundary at time $\dt$, $\n$ is the 
inward normal vector at $\x_0$, $\kappa$ is
the curvature, and $v$ is the normal velocity. 
Here, $E_{3/2}$ denotes the exponential integral of order 3/2:
\begin{equation}
E_{3/2}(x)=\int_1^{\infty}\frac{e^{-xt}}{t^{3/2}}dt \,. 
\end{equation}
\end{lemma}

\noindent {\em Proof:}
Beginning with (\ref{eq:simplifiedsloc}), we assume that the functions
$y(s,\tau)$ and $\sigma(s,\tau)$ have the following Taylor expansions:
\begin{eqnarray}\label{yexp}  
y(s,\tau) &=& \sum_{n,m} \alpha_{n,m} s^n(\dt-\tau)^m=\sum_{n,m}\frac{\alpha_{n,m}}{4^m}u^n z^{2m+n} \\
\label{sigmaexp}
\sigma(s,\tau) &=& \sum_{n,m} \beta_{n,m} s^n(\dt-\tau)^m=\sum_{n,m}\frac{\beta_{n,m}}{4^m}u^n z^{2m+n}
\end{eqnarray}
Note that 
\begin{align*}
\alpha_{0,0}&=y(0,\dt)=0 \\
\alpha_{1,0}&=y_s(0,\dt)=0 \\
\alpha_{0,1}&=-y_{\tau}(0,\dt)=-v \\
\alpha_{2,0}&=\frac{y_{ss}(0,\dt)}{2}=\frac{\kappa}{2} \, .
\end{align*}
The result now follows by 
substituting (\ref{yexp}) and (\ref{sigmaexp}) into 
$e^{-(r-y)^2/z^2}$ and $\sqrt{1+y_s^2}$ in equation (\ref{eq:simplifiedsloc}), 
using the facts that
\begin{align*}
\int_{-\infty}^{\infty} &u^{2n-1}e^{-u^2} du = 0 \\
\int_{-\infty}^{\infty} &u^{2n}e^{-u^2} dr = 
\sqrt{\pi}\cdot\frac{(2n-1)!!}{2^n} 
\end{align*}
$\qed$

\begin{remark}
This calculation can be carried out to higher order in $\dt$. 
However, the expression will involve
higher order derivatives of $y(s,t)$ and $\sigma$ and are not likely to 
be practical. The expansions above are sufficient for our purposes.
\end{remark}

\begin{remark}
The exponential integral of order $s$ is defined by
\[ E_s(x)=\int_1^{\infty}\frac{e^{-xt}}{t^s} dt. \]
It satisfy the recurrence: 
\[ E_{s+1}(x)=\frac{1}{s}(e^{-x}-xE_s(x)) \]
for $s>0$, and 
\[ E_s(x)=x^{s-1}\cdot \Gamma(1-s,x) \]
for $0<s<1$. Here $\Gamma(1-s,x)$
is the incomplete Gamma function of order $1-s$. 
In particular, we have
\[ E_{3/2}(x)=2(e^{-x}-\sqrt{x}\Gamma(\frac{1}{2},x)).
\]
This provides one method for computing $E_{3/2}(x)$, since 
incomplete Gamma functions are well studied special functions and there
is widely available software for their evaluation. 
\end{remark}

The asymptotic expansion for the double layer potential $D$ is
given by the following lemma.

\begin{lemma}\label{lemma:dlocasym}
Let $\Gamma(\tau)$ and $\mu(\x,\tau)$ be four times differentiable, then
\begin{equation}\label{eq:dlocasym}
\begin{split}
D[\mu](\x,\dt)=&-\sqrt{\frac{\dt}{\pi}}E_{3/2}\left(\frac{c^2}{4}\right)\frac{\kappa+v}{4}\mu(\x_0,\dt) \\
&-\frac{\sgn(c)}{2}\erfc\left(\frac{|c|}{2}\right)(1+c\sqrt{\dt}\cdot\frac{\kappa-v}{2})\mu(\x_0,\dt)+O(\dt),
\end{split}
\end{equation}
where $\erfc(x)=\frac{2}{\sqrt{\pi}}\int_x^{\infty} e^{-t^2} dt$ is the 
complementary error function. $\x$, $\x_0$, $\kappa$, $v$ and $c$ have
the same definition as in Lemma \ref{lemma:slocasym}.
\end{lemma}

\noindent {\em Proof:}
Using the same notation and change of variables as above, 
$D$ assumes the form
\begin{equation}\label{eq:dlocsimplified}
D[\mu](\x,\dt)=\frac{1}{\pi}\int_0^{2\sqrt{\dt}}
\int_{-\infty}^{\infty} e^{-u^2}e^{-\frac{(r-y)^2}{z^2}}
\frac{y-r-sy_s}{z^2}\mu(u,z) \, du \, dz.
\end{equation}
Assume now that $\mu(s,\tau)$ has the following Taylor expansion:
\begin{equation} \label{muexp}
\mu(s,\tau)=\sum_{n,m} \gamma_{n,m} s^n(\dt-\tau)^m=
\sum_{n,m}\frac{\gamma_{n,m}}{4^m} u^n z^{2m+n}
\end{equation}
Substituting (\ref{yexp}) and (\ref{muexp}) into $e^{-\frac{(r-y)^2}{z^2}}$ 
and $\frac{y-r-sy_s}{z^2}$ in equation (\ref{eq:dlocsimplified})
yields the desired result, using the fact that
\begin{equation}
\int_0^{2\sqrt{\dt}} \frac{r^2/z^2}{z^2} dz=\frac{\sqrt{\pi}}{2|r|}\erfc\left(\frac{|r|}{2\sqrt{\dt}}\right) \, .
\end{equation}
$\qed$

\begin{remark}
There is no essential difficulty in carrying out this calculation 
to higher order. In the case of $D[\mu]$, however, the expression is quite
involved, even to achieve an error of the order $O(\dt^{3/2})$.
It involves the derivatives $\mu_{\tau}$ and $\mu_{ss}$:
\begin{equation}\label{eq:dlocasym32}
\begin{split}
D[\mu](\x,\dt)&=\sqrt{\frac{\dt}{\pi}}E_{3/2}\left(\frac{c^2}{4}\right)\left(-\frac{\kappa+v}{4}+\frac{v^2-3\kappa^2}{8}c\sqrt{\dt}\right)\mu(\x_0,\dt) \\
&-\frac{\sgn(c)}{2}\erfc\left(\frac{|c|}{2}\right)\left(1+c\sqrt{\dt}\frac{\kappa-v}{2}+c^2\dt \frac{v^2+3\kappa^2-2v\kappa}{8}\right)\mu(\x_0,\dt) \\
&-\frac{c\dt}{\sqrt{\pi}}E_{3/2}\left(\frac{c^2}{4}\right)\left(\frac{v^2+3\kappa^2-2v\kappa}{16}\mu(\x_0,\dt)+\frac{\mu_{ss}(\x_0,\dt)-\mu_{\tau}(\x_0,\dt)}{4}\right) \\
& +O(\dt^{3/2}) \, .
\end{split}
\end{equation}
\end{remark}

\subsection{Asymptotics on the boundary\label{sec:two}} 

\vspace{.01in}

In the asymptotic formulae (\ref{eq:slocasym}) and 
(\ref{eq:dlocasym}), we may let $\x\rightarrow \x_0$ 
(i.e. $c\rightarrow 0$) in order to recover the
known asymptotics \cite{greengard_strain2,li09,lin_thesis} 
for points on the boundary $\Gamma(\dt)$ itself.
 For the sake of completeness, we include them here.

\begin{cor}\label{cor:slocasym0}
Let $\Gamma(\tau)$ and $\sigma(\x,\tau)$ be four times differentiable, and 
let $\x_0\in\Gamma(\dt)$ be a point on the boundary at time $\dt$. Then
\begin{equation}\label{eq:slocasym0}
S[\sigma](\x_0,\dt)=\sqrt{\frac{\dt}{\pi}}\sigma(\x_0,\dt)+O(\dt^{3/2}).
\end{equation}

\end{cor}

\noindent {\em Proof:}
In equation (\ref{eq:slocasym}), let $\x\rightarrow \x_0$ 
(i.e. $c=\frac{|\x-\x_0|}{\sqrt{\dt}}\rightarrow 0$).
The result follows from the facts that 
$\lim_{c\rightarrow 0}E_{3/2}(c)=E_{3/2}(0)=2$ and 
that the single layer is only weakly singular, so that
$$\lim_{x\rightarrow \x_0\in\Gamma(\dt)} S[\sigma](\x,\dt)=
S[\sigma](\x_0,\dt). $$
$\qed$

\begin{cor}\label{cor:dlocasym0}
Let $\Gamma(\tau)$ and $\mu(\x,\tau)$ be four times differentiable, and
let $\x_0\in\Gamma(\dt)$ be a point on the boundary at time $\dt$. Then
\begin{equation}\label{eq:dlocasym0}
D^{*}[\mu](\x_0,\dt)=-\sqrt{\frac{\dt}{\pi}}
\frac{\kappa+v}{2}\mu(\x_0,\dt)+O(\dt^{3/2}),
\end{equation}
where $D^{*}[\mu](\x_0,\dt)$ is defined in 
\eqref{eq:dl}.
\end{cor}

\noindent {\em Proof:}
The result follows from
equation (\ref{eq:dlocasym32}), letting $\x\rightarrow \x_0$ 
(i.e. $c=\frac{|\x-\x_0|}{\sqrt{\dt}}\rightarrow 0$). 
We also make use of the jump relations
$$\lim_{\x\rightarrow \x_0\in\Gamma(\dt), \x\in \Omega(\dt)} 
D[\sigma](\x,\dt)=-\frac{1}{2}\mu(\x_0,\dt)+D^{*}[\mu](\x_0,\dt)$$
$$\lim_{\x\rightarrow \x_0\in\Gamma(\dt), \x\in \Omega^c(\dt)} 
D[\sigma](\x,\dt)=\frac{1}{2}\mu(\x_0,\dt)+D^{*}[\mu](\x_0,\dt),$$
with the convention that when $\x$ is in the interior of the domain
($\x\in \Omega(\dt)$), we have $\sgn(c)=1$. When $\x\in \Omega^c(\dt)$ 
we have $\sgn(c)=-1$.
$\qed$

\begin{remark}
Note that, on the boundary, the order of accuracy for $D$ 
is already of the order $O(\dt^{3/2})$, including only the 
leading order term.
This is because all of the terms of order $O(\dt)$
in equation \ref{eq:dlocasym32} vanish as $c\rightarrow 0$.  
\end{remark}

\subsection{Generalization}
Sometimes it is convenient to have asymptotic formulas for the 
layer potentials that are
one step removed from the current time, which we will now denote by 
$2\dt$. That is, we seek approximations of
\begin{equation}
 S_{B}[\sigma](\x,2\dt)=\int_{0}^{\dt} \int_{\Gamma(\tau)} 
G(\x-\y,2\dt-\tau)\sigma(\y,\tau) ds_{\y}d\tau 
\end{equation}
\begin{equation}
 D_{B}[\mu](\x,2\dt)=\int_{0}^{\dt} \int_{\Gamma(\tau)} 
\frac{\partial}{\partial n_{\y}}G(\x-\y,2\dt-\tau)\mu(\y,\tau) ds_{\y}d\tau
\end{equation}

All of the transformations used above still apply and
we give the results here without repeating the derivation.

\begin{lemma}\label{lemma:sbtasym}
Let $\Gamma(\tau)$ and $\sigma(\x,\tau)$ be four times differentiable, then
\begin{equation}\label{eq:sbtasym}
S_B[\sigma](\x,2\dt)=\frac{1}{2}
\sqrt{\frac{\dt}{\pi}}
\left(\sqrt{2}E_{3/2}\left(\frac{c^2}{8}\right)-E_{3/2}\left(\frac{c^2}{4}\right) \right)
\left(1+\frac{\kappa-v}{2}\cdot c\sqrt{\dt}\right)
\cdot\sigma(\x_0,2\dt)+O(\dt^{3/2}),
\end{equation}
following the notation of Lemma (\ref{lemma:slocasym}).
\end{lemma}

\begin{lemma}\label{lemma:dbtasym}
Let $\Gamma(\tau)$ and $\mu(\x,\tau)$ be four times differentiable, then
\begin{equation}\label{eq:dbtasym}
\begin{split}
D_B[\mu](\x,2\dt)=&-\sqrt{\frac{\dt}{\pi}}\left(\sqrt{2}E_{3/2}\left(\frac{c^2}{8}\right)-E_{3/2}\left(\frac{c^2}{4}\right)\right)\frac{\kappa+v}{4}\mu(\x_0,2\dt) \\
&-\frac{\sgn(c)}{2}\left(\erfc\left(\frac{\sqrt{2}|c|}{4}\right)-\erfc\left(\frac{|c|}{2}\right)\right)(1+c\sqrt{\dt}\cdot\frac{\kappa-v}{2})\mu(\x_0,2\dt) \\
&+O(\dt),
\end{split}
\end{equation}
following the notations of lemma (\ref{lemma:dlocasym}).
\end{lemma}

For targets $\x_0\in\Gamma(2\dt)$ on the boundary, 
simply let $x\rightarrow 0$ in the above lemmas.

\begin{cor}\label{cor:sbtasym0}
Let $\Gamma(\tau)$ and $\sigma(\x,\tau)$ be four times differentiable, and
$\x_0\in\Gamma(2\dt)$ be a point on the boundary at time $\dt$, then:
\begin{equation}\label{eq:sbtasym0}
S_B[\sigma](\x_0,2\dt)=\frac{\sqrt{2}-1}{2}\sqrt{\frac{\dt}{\pi}}\sigma(\x_0,2\dt)+O(\dt^{3/2}) \, .
\end{equation}
\end{cor}

\begin{cor}\label{cor:dbtasym0}
Let $\Gamma(\tau)$ and $\mu(\x,\tau)$ be four times differentiable, and let
$\x_0\in\Gamma(2\dt)$ be a point on the boundary at time $2\dt$. Then
\begin{equation}\label{eq:dbtasym0}
D_B[\mu](\x_0,2\dt)=-\sqrt{\frac{\dt}{\pi}}(\sqrt{2}-1)\frac{\kappa+v}{2}\mu(\x_0,2\dt)+O(\dt^{3/2}) \, .
\end{equation}
\end{cor}

\section{Quadrature methods for local heat potentials} \label{sec_quad}
In the preceding section, we developed asymptotic methods for $S[\sigma]$ 
and $D[\mu]$.
Formally, these appear to be highly efficient,
but there are two major drawbacks to their use.
First, as we saw above, it is somewhat unwieldy to carry out 
asymptotic expansions to high order in $\dt$.
Moreover, those expressions involve impractical high order derivatives
of the density and the geometry.
A second problem is that the formal order of convergence is not 
manifested until $\dt$ is rather small, on the order of 
$\Delta x^2$, where $\Delta x$ is the spacing in the discretization of the 
boundary. This lack of resolution for larger values of $\dt$
is referred to as ``geometrically induced stiffness,"
and has been studied in detail in \cite{li09}. 

Leaving asymptotic methods aside for the moment, direct numerical quadrature
approximations require some care because of the singularity in time
introduced by the heat kernel.
To clarify the nature of the singularity in time, 
let us rewrite the local layer potentials in the following form:
\begin{equation} \label{eq:slocprime} 
 S[\sigma](\x,\dt)=\int_{0}^\dt 
\frac{1}{\sqrt{4\pi(\dt-\tau)}} B_S[\sigma](\x,\dt,\tau) d\tau,
\end{equation}
where
\begin{equation} \label{eq:bs}
B_S[\sigma](\x,\dt,\tau)=\int_{\Gamma(\tau)} \frac{e^{-\frac{|\x-\y|^2}{4(\dt-\tau)}}}{\sqrt{4\pi(\dt-\tau)}}\sigma(\y,\tau) ds_{\y},
\end{equation}
and 
\begin{equation} \label{eq:dlocprime}
D[\mu](\x,\dt)=\int_{0}^\dt \frac{1}{\sqrt{4\pi(\dt-\tau)}} B_D[\mu](\x,\dt,\tau) d\tau,
\end{equation}
where
\begin{equation}\label{eq:bd}
B_D[\mu](\x,\dt,\tau)=\int_{\Gamma(\tau)} \frac{e^{-\frac{|\x-\y|^2}{4(\dt-\tau)}}}{4\sqrt{\pi}(\dt-\tau)^{3/2}}(\x-\y)\cdot n_{\y}\mu(\y,\tau) ds_{\y}.
\end{equation}

A straightforward calculation shows that both $B_S[\sigma](\x,\dt,\tau)$ and 
$B_D[\mu](\x,\dt,\tau)$ are bounded as $\tau\rightarrow \dt$ 
and of the same order of smoothness as $\sigma$ and $\Gamma(\tau)$. 
Thus, the singularity in 
\eqref{eq:slocprime} and \eqref{eq:dlocprime} is indeed of the form
$\frac{1}{\sqrt{\dt-\tau}}$, as written.
There are an abundance of quadrature rules available to deal with 
integrands of this form.

\subsection{Partial product integration}
Because of the inverse square root singularity in 
\eqref{eq:slocprime},
a natural choice is product integration, which takes the form:
\begin{equation}
S[\sigma](\x,\dt)=\sqrt{\frac{\dt}{4\pi}}\sum_{j=0}^k w_j B_S[\sigma](\x,\dt,v_j)+e_{\dt}(k).
\end{equation}
Here, the nodes $v_j=\frac{j}{k}\dt$ are equispaced,
and the weights $w_0,\cdots, w_k$ are chosen so that
\begin{equation}
\int_{0}^\dt \frac{1}{\sqrt{\dt-\tau}}g(\tau) d\tau =\sqrt{\dt}\sum_{j=0}^k w_jg(\dt-v_j)
\end{equation}
is exact for $g(\tau)$ a polynomial of degree $\leq k$. It is straightforward to verify that the error $e_{\dt}(k)$
is of the order $k+3/2$. More precisely, we have: 
\begin{equation} \label{errPI}
e_{\dt}(k)\leq C\dt^{k+3/2} \|g^{(k+1)} \|_{\infty},
\end{equation}
for some constant $C$. Variants of this method using nonequispaced nodes 
include Gauss-Jacobi quadrature \cite{DR},
and hybrid Gauss-trapezoidal rules \cite{alpert_quad}. 
With a slight abuse of language, we refer to all such methods 
as {\em partial product integration} methods, 
since they take into account only part of the structure of the
heat kernel. An important feature of these methods is that each spatial 
integral is simply a convolution with a Gaussian and can be computed in 
linear time using the fast Gauss transform 
\cite{greengard_strain1,pfgt,tausch-fgt,fgtvol}.

Unfortunately, these methods also suffer from 
{geometrically induced stiffness}, as shown in 
\cite{li09}.
The mode of failure is that, even when the geometry and density are 
sufficiently well-resolved with a spatial grid with spacing $\Delta x$,
the formal order of accuracy of partial product integration is sometimes
not manifested until $\dt \approx \Delta x^2$.
While first observed as a geometric phenomenon near regions of high curvature,
it occurs even on straight boundaries with sharply peaked densities,
so that the term {\em geometric} is perhaps unfortunate. 

\subsection{Geometrically induced stiffness}

As a simple example, let us consider the 
computation of $S[\sigma]$, where the boundary $\Gamma$ is
chosen to be the x-axis, and the density $\sigma$ is chosen to be 
a (time-independent) Gaussian centered at the origin:
\begin{equation} \label{densbump}
\sigma(\y,\tau)=\sigma_d(y)=\frac{1}{\sqrt{4\pi d}}e^{-\frac{y^2}{4d}},\;\;\;(\y=(y,0)).
\end{equation}
Note that $\sigma(\y,\tau)$ is well-resolved with a spatial mesh whose spacing
is $\Delta x \approx \sqrt{d}$.
We choose the evaluation point to be $\x=(x,0)$. In this case, 
letting $t = \Delta t - \tau$,
$S[\sigma]$ takes the form:
\begin{equation} \label{eq:slocbump}
\begin{split}
S[\sigma](\x,\dt)&=\int_0^{\dt} \int_{-\infty}^{\infty} \frac{1}{4\pi t}e^{-\frac{(x-y)^2}{4 t}}
\sigma_d(y)dy dt\\
& =\int_0^{\dt}\frac{1}{\sqrt{4\pi t}} \int_{-\infty}^{\infty} \frac{1}{\sqrt{4\pi t}}e^{-\frac{(x-y)^2}{4 t}}\sigma_d(y)dy d t \\
& =\int_0^{\dt} \frac{1}{\sqrt{4\pi t}} g_d(x,t) dt,
\end{split}
\end{equation}
where
\begin{equation} \label{eq:slocbumpdens}
g_d(x,t)=\frac{1}{\sqrt{4\pi(t+d)}}e^{-\frac{x^2}{4(t+d)}} .
\end{equation}

It is straightforward to verify that
\begin{equation} 
\frac{\partial^k}{\partial t^k}g_d(x,t)=
\frac{1}{\sqrt{4\pi}}\left(\frac{1}{2}\right)^{2k}\left(\frac{1}{t+d}\right)^{k+1/2}
h_{2k}\left(\frac{x}{\sqrt{2t+d}}\right),
\end{equation}
where $h_{2k}$ is the Hermite function of order $2k$, defined as $h_{2k}(t)=D^{2k} e^{-t^2}$. At $x=0$,
for example, the maximum value of $\frac{\partial^k}{\partial t^k}g_d(x,t)$ is achieved at $t=0$,
and thus we have:
\begin{equation} \label{derivbound}
\|\frac{\partial^k}{\partial t^k}g_d(0,\cdot)\|_{\infty}=
\frac{2^k(2k-1)!!}{\sqrt{4\pi d}}\left(\frac{1}{4d}\right)^k .
\end{equation}

Suppose now that we make use of product integration with $k+1$ nodes in time.
Combining (\ref{derivbound}) with (\ref{errPI}), we obtain:
\begin{equation} \label{errorboundPPI}
\begin{split}
e_{\dt}(k)&\leq C\dt^{k+3/2}\frac{2^{k+1}(2k+1)!!}{\sqrt{\pi}}\left(\frac{1}{4d}\right)^{k+3/2}\\
&=C(k)\left(\frac{\dt}{4d}\right)^{k+3/2},
\end{split}
\end{equation}
where $C(k)$ is a constant that is independent of $\dt$ and $d$ 
but grows rapidly with $k$.
Thus, despite the fact that $\Delta x \approx \sqrt{d}$ is sufficient to 
resolve the density, we need $\dt \approx d \approx \Delta x^2$ in order
for the time quadrature to be accurate, 
so that $\dt\leq 4d$ in (\ref{errorboundPPI}). This is not a formal result:
partial product integration with larger time steps is, indeed, inaccurate 
in this regime.

\subsection{Adaptive Gaussian Quadrature}

A powerful method for handling end-point singularities is 
adaptive Gaussian quadrature. The basic idea for an integrable
singularity at $\tau = \dt$ is to subdivide the time
interval $[0,\dt]$ dyadically, so that each subinterval is separated from
the singular endpoint by its own length. For a precision 
$\epsilon$, it is easy to see that cutting off the last 
interval $[\dt-\epsilon^2,\dt]$ 
yields a truncation error of size $\epsilon$.
This follows from the fact that the norm of the operator 
$S[\sigma](\x,\epsilon)$ is of 
the order $O(\sqrt{\epsilon})$.
Dydadic refinement to this scale clearly requires approximately 
$\log_2 \epsilon^2$ subintervals.
In the elliptic setting, say for corner singularities using two-dimensional
boundary integral methods, it has been shown that
using a simple $n$th order Gauss-Legendre rule on each subinterval
yields ``spectral accuracy" with respect to $n$
\cite{BRS,triplepoint,HO,trefethen}.
That is, the error is of the order $O( e^{-n} \log_2 \epsilon)$. 
We now show that the same principle can be applied to heat potentials,
in a manner that is high order and that also overcomes geometrically 
induced stiffness.
For the sake of simplicity, we restrict our attention to the 
same model problem considered above.

\begin{theorem} \label{thm:bumpdens}
Let $S^{\epsilon}[\sigma](\x,\dt)$ be the single layer potential
truncated at $\dt - \epsilon^2$:
\begin{equation} \label{eq:sloceps}
S^{\epsilon}[\sigma](\x,\dt)=\int_{0}^{\dt-\epsilon^2} 
\int_{-\infty}^\infty \frac{1}{4\pi(\dt-\tau)}
e^{-\frac{(x-y)^2}{4(\dt-\tau)}}\sigma_d(y) \, dy d\tau,
\end{equation}
where $\sigma_d(y)$ is given by (\ref{densbump}), and the evaluation point
is chosen to be $\x=(x,0)$. Suppose now that the interval 
$[0,\dt-\epsilon^2]$ has a dyadic decomposition of the form:
\begin{equation} \label{eq:dyadic}
[0,\dt-\epsilon^2]=\cup_{i=0}^{N-1} [\dt-2^i\epsilon^2,\dt- 2^{i+1}\epsilon^2],
\end{equation}
where $\dt=2^N\epsilon^2$, 
and that on each dyadic subinterval,
the Gauss Legendre rule with $n$ quadrature nodes is applied. 
Then the error $E_n$ in applying composite Gauss-Legendre quadrature
satisfies the error bound:
\begin{equation} 
E_n\leq C\log_2\left(\frac{\dt}{\epsilon^2}\right)\frac{n^{1/4}}{16^n}.
\label{errorestgauss}
\end{equation}
\end{theorem}

\noindent {\em Proof:}
Letting $t = \dt-\tau$,
from (\ref{eq:slocbump}) and (\ref{eq:slocbumpdens}), 
we have that
\begin{equation} \label{eq:slocbump2}
S^{\epsilon}[\sigma](\x,\dt)=\int_{\epsilon^2}^{\dt} \frac{1}{\sqrt{4\pi t}} g_d(x,t) dt,
\end{equation}
where
\begin{equation} \label{eq:slocbumpdens2}
\begin{split}
g_d(x,t)&=\frac{1}{\sqrt{4\pi(t+d)}}e^{-\frac{x^2}{4(t+d)}} \\
&=\frac{1}{\sqrt{4\pi(t+d)}}h_0\left(\frac{x}{2\sqrt{t+d}}\right),
\end{split}
\end{equation}
and $h_0$ is the Hermite function of order zero.
For convenience of notation, we let
\begin{equation}
F_d(x,t)=\frac{1}{\sqrt{4\pi t}}g_d(x,t) \, .
\end{equation}

We now rewrite the integral as a sum of integrals on dyadic subintervals:
\begin{equation}
S^{\epsilon}[\sigma](\x,\dt)=\sum_{i=1}^{N-1} \int_{2^i\epsilon^2}^{2^{i+1}\epsilon^2} F_d(x,t) dt,
\end{equation}
and denote by $E_n^{i}$ the quadrature error on the $i$th subinterval.
Thus,
\begin{equation}
E_n^i= 
\left\| \int_{2^i\epsilon^2}^{2^{i+1}\epsilon^2} F_d(\cdot,t) dt
-\sum_{j=1}^n w_{ij} F_d(\cdot,t_{ij}) \right\|_{\infty},
\end{equation}
where $\{t_{ij}\}_{j=1}^n$ and $\{w_{ij}\}_{j=1}^n$ 
are the Legendre nodes and weights on the $i$th
subinterval. In order to apply the standard estimate 
for Gauss Legendre quadrature (Lemma \ref{lemma:GLQuad} in the appendix),
we need an upper bound for the $(2n)^{th}$ derivative of $F_d(x,t)$ with
respect to time.

Since $g_d(x,t)$ satisfies the one dimensional heat equation on 
$(-\infty,\infty)\times (-d, \infty)$,
we have:
\begin{equation}
D_{t} g_d(x,t)=D_x^2 g_d(x,t), \;\;(\forall (x,t)\in(-\infty,\infty)\times (0, \infty)),
\end{equation}
so that
\begin{equation}
\begin{split}
&D_{t}^k g_d(x,t)=D_x^{2k} g_d(x,t) \\
&=\frac{1}{\sqrt{4\pi}}\left(\frac{1}{2}\right)^{2k}\left(\frac{1}{t+d}\right)^{k+1/2}
h_{2k}\left(\frac{x}{2\sqrt{t+d}}\right). 
\end{split}
\end{equation}

Applying Cramer's inequality (Lemma \ref{lemma:Cramer}),
we obtain:
\begin{equation}
|D_{t}^k g_d(x,t)|\leq C\left(\frac{1}{2}\right)^k \left(\frac{1}{t}\right)^{k+1/2}
\sqrt{(2k)!}  \, .
\end{equation}

From Corollary \ref{cor:factorial} (see the appendix), this implies
\begin{equation}
|D_{t}^k g_d(x,t)|\leq \tilde{C} n^{1/4}
\left|D_{t}^k \left(\frac{1}{\sqrt{t}}\right) \right|,
\end{equation}
which means that up to a mild growth factor in $n$, 
the $k$th derivative of $g_d(x,t)$ is
controlled by that of $\frac{1}{\sqrt{t}}$.
Returning to the integrand in (\ref{eq:slocbump2}),
Leibniz's rule leads to 
\begin{equation}
\begin{split}
|D_{t}^{2n} F_d(x,t)|&=\sum_{k=0}^{2n} {2n \choose k} 
\left|D_{t}^{2n-k}
\left(\frac{1}{\sqrt{4\pi t}}\right) \right|\cdot|D_{t}^k g_d(x,t)| \\
& \leq Cn^{1/4}\sum_{k=0}^{2n} {2n\choose k} 
\left|D_{t}^{2n-k}\left(\frac{1}{\sqrt{t}}\right) \right|
\cdot \left|D_{t}^{k}\left(\frac{1}{\sqrt{t}}\right) \right| \\
& = Cn^{1/4} \left|D_{t}^{2n}\left(\frac{1}{t}\right) \right|,
\;\;\;(\forall t>0,\;\;\forall x\in\bR).
\end{split}
\end{equation}
(The last step follows from the Leibnitz rule applied to 
$\frac{1}{t} = \frac{1}{\sqrt{t}} \cdot \frac{1}{\sqrt{t}}$.)

Restricting $t$ to the $i$th subinterval 
$[2^i\epsilon^2, 2^{i+1}\epsilon^2]$ we obtain:
\begin{equation}
D_{t}^{2n} F_d(x,t) \leq Cn^{1/4}(2n)!\left(\frac{1}{2^i\epsilon^2}\right)^{2n+1}
\;\;\;\forall x\in\bR.
\end{equation}
Combined with the standard error estimate for Gauss-Legendre quadrature,
this yields the following bound:
\begin{equation}
E_n^i\leq Cn^{1/4}\left(\frac{1}{16}\right)^n.
\end{equation}

It remains only to sum up $E_n^i$ over all subintervals to obtain 
\begin{equation}
\label{allerrs}
E_n \leq\sum_{i=1}^{N-1} E_n^i
=C\log_2\left(\frac{\dt}{\epsilon^2}\right)\frac{n^{1/4}}{16^n}.
\end{equation}
$\qed$

The above theorem shows the exponential convergence rate 
of adaptive Gaussian quadrature
for a specific choice of density function. It can easily be generalized 
to an arbitary smooth density
function in the following manner. Suppose that $\sigma(\y,t)$ is 
time dependent. We then expand 
the density function as a Taylor series in time around $t=\dt$:
\begin{equation}
\begin{split}
\sigma(\y,\tau)&=\sigma_0(\y)+\sigma_1(\y)(\dt-\tau)+\frac{1}{2}\sigma_2(\y)(\dt-\tau)^2 \\
&+\cdots+\frac{1}{(k-1)!}\sigma_{k-1}(\y)(\dt-\tau)^{k-1}+O((\dt-\tau)^k).
\end{split}
\end{equation}
Denoting by $\tilde \sigma_k(\y,\tau)$ the first k terms in the expansion, we have
\begin{equation}
S[\sigma](\x,\dt)=S[\tilde \sigma_k](\x,\dt)+O(\dt^{k+1/2}).
\end{equation}
A similar analysis to the estimate
\eqref{errorestgauss} in  the proof of Theorem \ref{thm:bumpdens}
leads to a total error of the form:
\begin{equation} 
E_{n,k}\leq C_1\dt^{k+1/2}+C_2 \log_2 \left( \frac{\Delta t}{\epsilon^2} \right) \,
\frac{n^{1/4}}{16^n} + O(\epsilon).
\label{etotal}
\end{equation}
The first term on the right-hand side is due to the 
$k$th order accurate approximation of the density as a function of time.
The second term is due to the error in adaptive Gaussian 
quadrature, as in \eqref{allerrs}, 
and the last term is due to ignoring the contribution of the 
last time interval $[\Delta t - \epsilon^2,\Delta t]$ to the 
single layer potential.

When the boundaries are 
curves, the calculation is more complicated, involving time derivatives of 
a spatial convolution on a possibly nonstationary domain.
Instead of entering into a detailed analysis, we illustrate the 
convergence behavior via numerical examples in the next section. 
Finally, we note that the double layer potential can be treated in 
essentially the same manner.

\subsection{A hybrid asymptotic/quadrature-based method}

In our initial analysis above, we integrated in time over the interval
$[\dt - \epsilon^2,\dt]$ with a truncation error of the order $O(\epsilon)$.
Let us, however, reconsider the use of the asymptotic approach, not on the 
entire interval $[0,\dt]$ but simply on 
$[\dt-\delta,\dt]$. For small $\delta$, geometrically induced stiffness
is no longer an issue and the cost of evaluation is negligible - requiring
only one kernel evaluation for each target point.
Thus, for the single layer potential, we decompose it in the form:
\begin{equation} \label{eq:slocsplit}
\begin{split}
 S[\sigma](\x,\dt)&=\int_{\dt-\delta}^\dt \int_{\Gamma(\tau)} 
G(\x-\y,\dt-\tau)\sigma(\y,\tau) ds_{\y}d\tau\\
 &+\int_{0}^{\dt-\delta} \int_{\Gamma(\tau)} 
G(\x-\y,\dt-\tau)\sigma(\y,\tau) ds_{\y}d\tau \\
 &:=S^{(1)}[\sigma,\delta](\x,\dt)+S^{\sqrt{\delta}}[\sigma](\x,\dt),
\end{split}
\end{equation}
where $S^{\sqrt{\delta}}[\sigma]$ is defined in 
\eqref{eq:sloceps}.
For a user-specified tolerance of $\epsilon$,
we choose $\delta$ so that the asymptotic formula (\ref{eq:slocasym}) 
for $S^{(1)}$ is accurate to a tolerance of $\epsilon$.
(This requires that $\delta$ be of the order $\Delta x^2$ or smaller
to avoid geometrically induced stiffness.)
$S^{\sqrt{\delta}}$ can then be treated
as above with {\em fewer} subintervals than before, since dyadic subdivision
is needed only until the smallest interval is of length $\delta$.

\subsection{Local quadrature with an exponentially graded mesh}

While we have shown that 
adaptive Gaussian quadrature on dyadic intervals
is robust and overcomes the problem of geometrically induced stiffness,
the total number of quadrature nodes in time is still relatively large.
We now seek a method for reducing the cost by designing a ``single 
panel" quadrature rule that maintains the exponential clustering toward the 
singularity inherent in the dyadic subdivision process. 
In the elliptic setting (for domains with corner singularities), this 
is discussed in \cite{ATK,BRS,HO}.
A simple mapping with this effect is obtain by setting
$e^{-u} = t-\tau$.  Under this change of variables,
$S^{\sqrt{\delta}}[\sigma](\x,\dt)$ becomes 
\begin{equation} \label{sloc2new}
S^{\sqrt{\delta}}[\sigma](\x,\dt)=
\frac{1}{4\pi}\int_{-\log{\dt}}^{-\log{\delta}}\int_{\Gamma(t-e^{-u})} e^{-\frac{|\x-\y|^2}{4e^{-u}}} \sigma(\y,t-e^{-u}) ds_{\y} du
\end{equation}
which is a smooth integral on the interval $[-\log{\dt}, -\log{\delta}]$. 
Applying Gauss-Legendre quadrature with $n$ nodes, we have
\begin{equation}
S^{\sqrt{\delta}}[\sigma](\x,\dt)= 
\frac{1}{4\pi} \sum_{j=1}^n \omega_j \int_{\Gamma(\dt-e^{-u_j})} 
e^{-\frac{|\x-\y|^2}{4e^{-u_j}}} \sigma(\y,\dt-e^{-u_j}) ds_{\y} 
+ e_S(\delta,\dt,n)
\end{equation}
where $\{u_j\}_{j=1}^n$ and $\{\omega_j\}_{j=1}^n$ are Legendre nodes and 
weights scaled to $[-\log{\dt}, -\log{\delta}]$. 
Here,  
$e_S(\delta,\dt,n)$ denotes the truncation error.

The double layer potential can be treated in a similar manner.
With the same
change of variables as above, namely $e^{-u} = t-\tau$, we obtain
\begin{equation}
D^{\sqrt{\delta}}[\mu](\x,\dt)=\frac{1}{8\pi}\int_{-\log{\dt}}^{-\log{\delta}}\int_{\Gamma(t-e^{-u})} e^{-\frac{|\x-\y|^2}{4e^{-u}}}
\mu(\y,t-e^{-u}) \frac{(\x-\y)\cdot n_{\y}}{e^{-u}} ds_{\y} du.
\end{equation}
Applying Gauss-Legendre quadrature with $n$ nodes yields
\begin{equation}
D^{\sqrt{\delta}}[\mu](\x,\dt)=
\frac{1}{8\pi}\sum_{j=1}^N \omega_j\int_{\Gamma(t-e^{-u_j})} e^{-\frac{|\x-\y|^2}{4e^{-u_j}}}
\mu(\y,t-e^{-u_j}) \frac{(\x-\y)\cdot n_{\y}}{e^{-u_j}} ds_{\y} + 
e_D(\delta,\dt,n).
\end{equation}
Here, $e_D(\delta,\dt,n)$ denotes the truncation error for the double layer 
potential.

To analyze the error for a single Gauss-Legendre panel with an 
exponentially graded mesh, 
let us consider as a model problem the calculation of
\begin{equation} \label{modelquad}
I(\delta,\dt)=\int_{\delta}^{\dt} \frac{1}{\sqrt{t}} \, dt.
\end{equation}
This has the same near singularity at
$t=0$ as $S^{\sqrt{\delta}}$ or $D^{\sqrt{\delta}}$ above.
If we carry out dyadic decomposition, as in (\ref{eq:dyadic}), 
and apply an $n$-point Gauss-Legendre rule on each dyadic interval,
we get the error bound:
\begin{equation} \label{eq:error_dyadic}
e_n(\delta,\dt)\leq C\cdot\frac{\sqrt{\dt}-\sqrt{\delta}}{\sqrt{n}\cdot16^n}.
\end{equation}
Note that in this case the total number of dyadic intervals is $\log_2(\frac{\dt}{\delta})$.
With $n$ quadrature nodes on each, 
it leads to $N=n\log_2(\frac{\dt}{\delta})$ nodes in total.

If we use instead the exponential change of variables, the
standard error estimate for Gauss-Legendre quadrature leads to:
\begin{equation}
I(\delta,\dt)=\sum_{j=1}^n \omega_j\cdot e^{-\tau_j/2} + e_n(\delta,\dt),
\end{equation}
where 
\begin{equation}\label{eq:error_exp}
e_n(\delta,\dt) \leq C\cdot\frac{\sqrt{\dt}}{n} \left(\frac{e\cdot\log(\dt/\delta)}{16\cdot n}\right)^{2n}.
\end{equation}
Both (\ref{eq:error_dyadic}) and (\ref{eq:error_exp}) suggest superalgebraic 
convergence in $n$ and low order accuracy in $\dt$.
This works very well in practice: the total error in the evaluation 
of a layer potential over a single time step is governed by the 
order of accuracy with which the density is approximated {\em and} 
the precision of the local rule determined by the parameter $n$ in 
\eqref{eq:error_exp}, as in \eqref{etotal} above.

For a concrete example, let us assume 
$\dt=10^{-2}$ and that the desired tolerance is $\epsilon=10^{-9}$.
We set $\delta = \epsilon$.
In this case we have $\log_2(\frac{\dt}{\delta})\approx 23$.
Adaptive Gaussian quadrature requires $n=9$ on each dyadic interval, 
meaning that the total number of quadrature
nodes is $n\log_2(\frac{\dt}{\delta})=207$, while the change of variables 
requires only $n=12$.
Rather than go through a detailed analysis of the full single and 
double layer potentials, numerical examples in the next section show that
the graded mesh performs as well as dyadic refinement with many 
fewer nodes, and that the errors
$e_S(\delta,\dt,n)$ and $e_D(\delta,\dt,n)$ decay superalgebraically with $n$.

\section{Numerical examples} \label{sec_num}
\begin{figure}[htbp] 
\centering
\includegraphics[width=.95\textwidth]{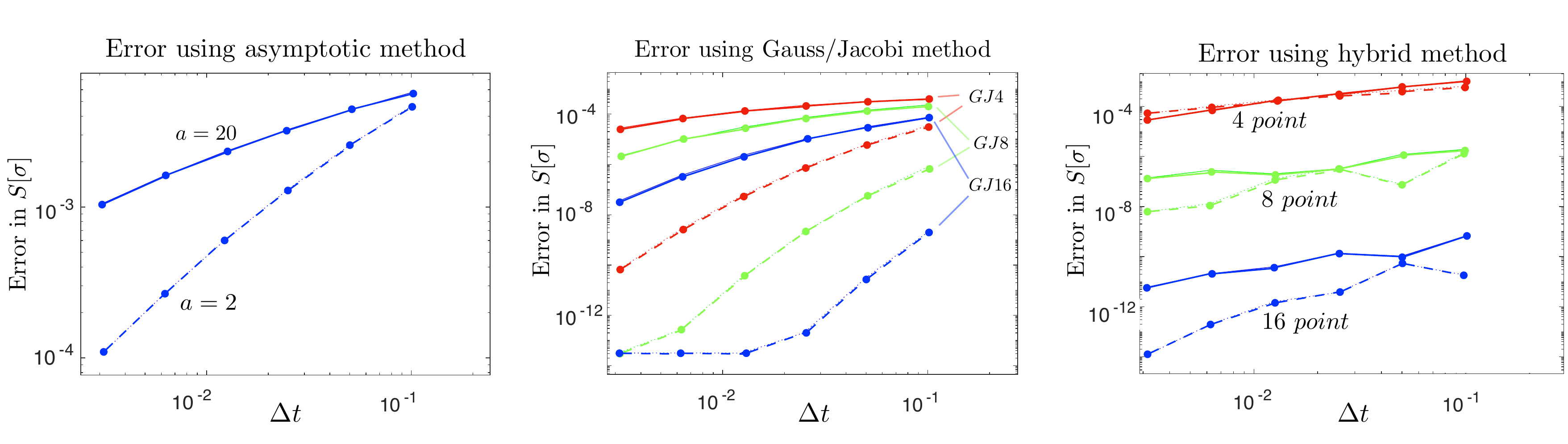} 
\caption{Comparison of quadrature methods for the single layer potential on 
a parabola with varying curvature.
The dashed lines correspond to the low curvature case ($a=2$) in all 
panels, while the solid lines correspond to the high curvature case ($a=20$).
}
\label{fig:exper1}
\end{figure}

For our first example, we consider as a boundary the parabola given by:
\begin{equation}
\begin{cases}
y_1(\lambda,t)= \lambda, \\
y_2(\lambda,t)=a\lambda^2,
\end{cases}
(-2\pi \leq \lambda \leq 2\pi)
\end{equation}
so that $2a$ is the curvature at the origin $\x_0 = (0,0)$.
We compute the single layer potential $S[\sigma](\x_0,\Delta t)$ with a 
constant density $\sigma(\y,t)\equiv 1$.
We carry out the experiment for $a=2$ and $a=20$, and use
as the time quadrature three different methods: 1) the asymptotic formula, 2) Gauss-Jacobi quadrature with four,
eight, and sixteen nodes, 3) our hybrid scheme with four, eight, and sixteen nodes. (We choose a small enough $\delta$
to give 13 digits of accuracy in the asymptotic regime.) 
All spatial integrals are
computed in this section to high precision, so that the errors in the examples come only from the time quadrature.
In Fig. \ref{fig:exper1}, quadrature errors are plotted for a wide range of $\Delta t$.

The results show that both the asymptotic formula and 
Gauss-Jacobi quadrature are very sensitive to the geometry.
For the case with low curvature, they perform relatively well, while for the case with high curvature, they are
inaccurate, especially at large values of $\Delta t$. Even for smaller values of $\Delta t$, they have still
not entered their rapidly convergent regime. 
The performance of the hybrid scheme is quite different. While low
order in $\Delta t$, it converges superalgebraically in the number of quadrature nodes, consistent with our analysis
above. With 16 quadrature nodes,
we are able to achieve 10 digits of accuracy for a wide range of $\Delta t$, even for the high curvature case,
where the other two methods perform poorly.

\begin{figure}[htbp] 
\centering
\includegraphics[width=.95\textwidth]{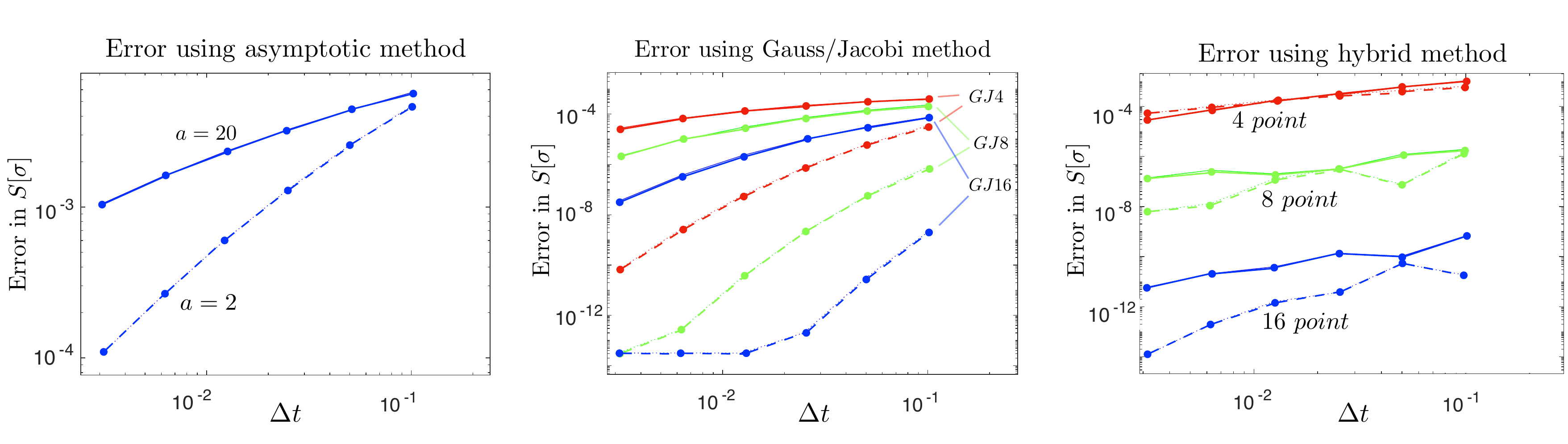} 
\caption{Comparison of quadrature methods for the single layer potential on 
a straight line segment with oscillatory density.
For all panels,
the dashed lines correspond to the case $k=10$ in \eqref{sigmadef},
while the solid lines correspond to $k=100$.
}
\label{fig:exper2}
\end{figure}

For our second example, we consider the computation of the single layer potential on a straight line segment $[-1,1]\times \{0\}$,
where the density is chosen to be an oscillatory function:
\begin{equation}
\label{sigmadef}
\sigma(\y,t)=\cos(2k\pi y_1).
\end{equation}
Errors for the same set of methods are plotted in Fig. \ref{fig:exper2}.
Again, the asymptotic formula and 
Gauss-Jacobi quadrature are sensitive to the frequency of the oscillation,
while the hybrid scheme is much more robust.

As a final example, we consider a more general task: computation 
of the double layer potential $D[\mu](\x,\Delta t)$
on a moving ellipse given by
\begin{equation}
\begin{cases}
y_1(\theta,t)= 20\cos(\theta)+1.5t, \\
y_2(\theta,t)= \sin(\theta),
\end{cases}
(0 \leq \theta \leq 2\pi)
\end{equation}
The density is chosen to be
\begin{equation}
\label{mudef}
\mu(\y,t)=\cos(y_1 t)+\sin(10t).
\end{equation}

\begin{figure}[htbp] 
\centering
\includegraphics[width=.95\textwidth]{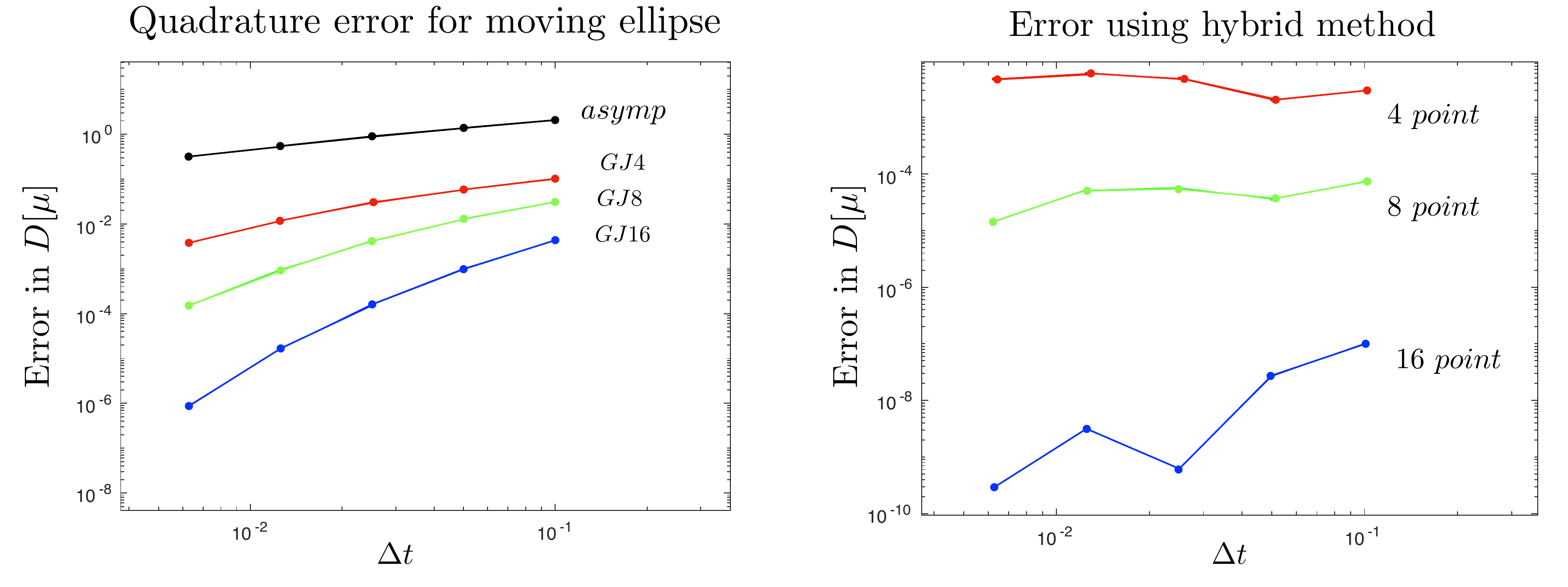} 
\caption{Comparison of quadrature methods for the double layer potential 
on a moving ellipse with oscillatory density given by \eqref{mudef}.
}
\label{fig:exper3}
\end{figure}

We evaluate the double layer potential at $\y=(21.5,0.0)$ and $t=1.0$ using the same set of methods. Errors are plotted in Fig. \ref{fig:exper3}.

\section{Conclusions} \label{sec_conclusions}

We have developed a 
new method for the evaluation of layer heat 
potentials in two dimensions. By making use of an exponential change of 
variables, we overcome the phenomenon of 
``geometrically-induced stiffness," which prevents the robust application
of high order Gauss-Jacobi type quadrature rules.
In our hybrid scheme, we combine a local asymptotic approximation with
Gauss-Legendre quadrature in the transformed time variable. The corresponding
spatial boundary integral operators involve only Gaussian kernels, permitting 
the application of the fast Gauss transform \cite{fgtvol}.
The scheme is easy to use with moving boundaries and precisely
the same rule can be used in three dimensions with respect to the time 
variable. A full solver for the heat equation in fixed and moving
geometries, incorporating the scheme of the present paper, will be 
described at a later date.

\section*{Acknowledgments}
We would like to thank Alex Barnett and Shidong Jiang for 
several useful conversations.

\appendix
\section{Stirling's formula, Cramer's inequality and Gauss-Legendre quadrature}

In the proof of Theorem \ref{thm:bumpdens}
we make use of Stirling's formula.

\begin{equation}
\label{lemma:Stirling}
\sqrt{2\pi}\left(\frac{n}{e}\right)^n\sqrt{n}
\leq n!
\leq e\left(\frac{n}{e}\right)^n\sqrt{n}.
\end{equation}

From this, it is staightforward to derive the following

\begin{cor} \label{cor:factorial}
Let $n\in\bN$, we have:
\begin{equation}
\sqrt{(2n)!} \leq C n^{1/4}(2n-1)!!\;\;,
\end{equation}
where $C>0$ is a constant.
\end{cor}

We also use Cramer's inequality \cite{hille}.

\begin{lemma}\label{lemma:Cramer}
Let $h_n(t)$ be the $n-th$ order Hermite function, defined by
\begin{equation}
h_n(t)=(-1)^n D^n e^{-t^2}.
\end{equation}
Then
\begin{equation}
|h_n(t)|\leq K 2^{n/2}\sqrt{n!}e^{-t^2/2},
\end{equation}
where $K$ is some constant with numerical value $K\leq 1.09$.
\end{lemma}

The following lemma is a direct consequence of Leibniz's 
product rule for differentiation.

\begin{lemma}
Let $F(t)=\frac{1}{t}$ and $f(t)=\frac{1}{\sqrt{t}}$. Then we have:
\begin{equation} \label{eq:leib}
F^{(2n)}(t)=\sum_{k=0}^{2n} {2n \choose k}|f^{(2n-k)}(t)|\cdot|f^{(k)}(t)|.
\end{equation}
\end{lemma}

Finally, we state the standard error estimate for Gauss-Legendre quadrature 
\cite{DR}.

\begin{lemma}\label{lemma:GLQuad}
Let $f\in C^{2n}([a,b])$ and
let $\{x_1,\cdots,x_n\}$ and $\{w_1,\cdots,w_n\}$ be the 
Gauss-Legendre nodes and weights scaled to 
$[a,b]$. If we denote the quadrature error by $E_n(f)$, we have:
\begin{equation} 
\begin{split}
E_n(f)&=\int_a^b f(x)dx-\sum_{k=1}^n w_k f(x_k) \\
&= \frac{(b-a)^{2n+1}(n!)^4}{(2n+1)[(2n)!]^3}f^{(2n)}(\xi),
\end{split}
\end{equation}
where $\xi\in[a,b]$.
\end{lemma}

\bibliographystyle{siam}
\bibliography{ref.bib}

\newcommand{\noopsort}[1]{} \newcommand{\printfirst}[2]{#1}
  \newcommand{\singleletter}[1]{#1} \newcommand{\switchargs}[2]{#2#1}
\begin{thebibliography}{10}

\bibitem{alpert_quad}
{\sc B.~K. Alpert}, {\em Hybrid {G}auss-trapezoidal quadrature rules}, SIAM J.
  Sci. Comput., 20 (1999), pp.~1551--1584 (electronic).

\bibitem{ATK}
{\sc K.~E. Atkinson}, {\em {The Numerical Solution of Integral Equations of the
  Second Kind}}, Cambridge University Press, New York, NY, 1997.

\bibitem{BRS}
{\sc J.~Bremer, V.~Rokhlin, and I.~Sammis}, {\em Universal quadratures for
  boundary integral equations on two-dimensional domains with corners}, J.
  Comput. Phys., 229 (2010), pp.~8259--8280.

\bibitem{DR}
{\sc P.~J. Davis and P.~Rabinowitz}, {\em Methods of numerical integration},
  Academic Press, San Diego, 1984.

\bibitem{triplepoint}
{\sc L.~Greengard and J.-Y. Lee}, {\em Stable and accurate integral equation
  methods for scattering problems with multiple interfaces in two dimensions},
  J. Comput. Phys., 231 (2012), pp.~2389--2395.

\bibitem{greengard_lin}
{\sc L.~Greengard and P.~Lin}, {\em Spectral approximation of the free-space
  heat kernel}, Appl. Comput. Harmon. Anal., 9 (2000), pp.~83--97.

\bibitem{greengard_strain2}
{\sc Leslie Greengard and John Strain}, {\em A fast algorithm for the
  evaluation of heat potentials}, Comm. Pure Appl. Math., 43 (1990),
  pp.~949--963.

\bibitem{greengard_strain1}
\leavevmode\vrule height 2pt depth -1.6pt width 23pt, {\em The fast {G}auss
  transform}, SIAM J. Sci. Statist. Comput., 12 (1991), pp.~79--94.

\bibitem{guentherlee}
{\sc R.~B. Guenther and J.~W. Lee}, {\em Partial differential equations of
  mathematical physics and integral equations}, Prentice Hall, Inglewood
  Cliffs, New Jersey, 1988.

\bibitem{HO}
{\sc J.~Helsing and R.~Ojala}, {\em Corner singularities for elliptic problems:
  integral equations, graded meshes, and compressed inverse preconditioning},
  J. Comput. Phys., 227 (2008), pp.~8820--8840.

\bibitem{hille}
{\sc E.~Hille}, {\em A class of reciprocal functions}, Annals of Mathematics:
  Second Series, 27 (1926), pp.~427--464.

\bibitem{li09}
{\sc Jing-Rebecca Li and Leslie Greengard}, {\em High order accurate methods
  for the evaluation of layer heat potentials}, SIAM Journal on Scientific
  Computing, 31 (2009), pp.~3847--3860.

\bibitem{lin_thesis}
{\sc P.~Lin}, {\em On the numerical solution of the heat equation in unbounded
  domains}, PhD thesis, New York University, 1993.

\bibitem{pogorzelski}
{\sc W.~Pogorzelski}, {\em Integral equations and their applications}, Pergamon
  Press, Oxford, 1966.

\bibitem{pfgt}
{\sc R.~S. Sampath, H.~Sundar, and S.~Veerapaneni}, {\em Parallel fast gauss
  transform}, in SC '10: Proceedings of the ACM/IEEE International Conference
  for High Performance Computing, Networking, Storage and Analysis, New
  Orleans, LA, 2010, pp.~1--10.

\bibitem{strain_adapheat}
{\sc J.~Strain}, {\em Fast adaptive methods for the free-space heat equation},
  SIAM J. Sci. Comput., 15 (1994), pp.~185--206.

\bibitem{tauschheat1}
{\sc J.~Tausch}, {\em A fast method for solving the heat equation by layer
  potentials}, J. Comput. Phys.,  (2007).

\bibitem{tausch-fgt}
{\sc J.~Tausch and A.~Weckiewicz}, {\em Multidimensional fast {G}auss
  transforms by {C}hebyshev expansions}, SIAM J. Sci. Comput., 31 (2009),
  pp.~3547--3565.

\bibitem{trefethen}
{\sc L.~N. Trefethen}, {\em Numerical computation of the schwarz-christoffel
  transformation}, SIAM J. Sci. Stat. Comput., 1 (1980), pp.~82--102.

\bibitem{fgtvol}
{\sc J.~Wang and L.~Greengard}, {\em An adaptive fast {G}auss transform in two
  dimensions}, SIAM J. Sci. Comput., to appear (arXiv:1712.00380) (2017).

\end{thebibliography}
\end{document}